\documentclass[twoside,12pt]{article}
\usepackage{amssymb,amsmath,bm,euscript}
\usepackage{graphicx,color,caption}

\usepackage[colorlinks=true]{hyperref}

\usepackage{algorithm,algpseudocode}
\usepackage{braket}

\newtheorem{proposition}{Proposition}[section]
\newtheorem{theorem}[proposition]{Theorem}
\newtheorem{lemma}[proposition]{Lemma}

\newcommand{\qed}{\hphantom{.}\hfill $\Box$\medbreak}

\def\x{{\bf x}}
\def\y{{\bf y}}
\def\z{{\bf z}}

\def\uu{{\bf u}}
\def\vv{{\bf v}}
\def\0{{\bf 0}}

\title{\bf{Generating Hypergraphs, Decomposability and Classification of Two-Step Nilpotent Lie Algebras}}%\thanks{This research  }
\author{ \hspace{1mm} Shenglong Hu\thanks{Department of Mathematics, Hangzhou Dianzi University, Hangzhou 310018 China; ({\tt shenglonghu@hdu.edu.cn}). This author's work was supported by NSFC (Grant No.  11771328).}, \ Liqun Qi\thanks{Department of Applied
    Mathematics, The Hong Kong Polytechnic University, Hung Hom,
    Kowloon, Hong Kong; ({\tt liqun.qi@polyu.edu.hk}).  This author's work was supported by the Hong Kong
    Research Grant Council (Grant No.  PolyU 15300715, 15301716 and 15300717).}
    \ and\
    Honglian Zhang\thanks{Department of Mathematics, Shanghai University, Shanghai, China; ({\tt hlzhangmath@shu.edu.cn}).  This author's work was supported by NSFC (Grant No.  11871325).}
    }

\begin{document}
\date{\today}
\maketitle

\begin{abstract}
In 1973, Gauger proposed a generator-relation method and a duality theory for two-step nilpotent Lie algebras.   Based upon these, he classified two-step nilpotent Lie algebras of dimension $8$.     In 1999, Galitski and Timashev continued this approach to classify two-step nilpotent Lie algebras of dimension $9$.   Their results were partially improved by Ren and Zhu in 2011, Yan and Deng in 2013.
 Some decomposable two-step nilpotent Lie algebras were excluded in the case of dimension $8$.      In this paper, we define generating hypergraph for a two-step nilpotent Lie algebra.  The two-step nilpotent Lie algebra is decomposable if and only if its generating hypergraph is not connected under certain bases.  Using this result, we identify some decomposable two-step nilpotent Lie algebras in dimension $9$.  We give a direct proof that the five two-step nilpotent Lie algebras for dimension $8$, classified by Ren and Zhu in 2011, are all indecomposable.   We also introduce a conventional nomenclature for two-step nilpotent Lie algebras of dimension $n = 8, 9$, classified by Ren, Zhu, Yan and Deng, etc.


%15A69:LINEAR AND MULTILINEAR ALGEBRA; MATRIX THEORY- Multilinear algebra, tensor products
%53A45:Classical differential geometry-Vector and tensor analysis
%47A05: Operator Theory-General (adjoints, conjugates, products, inverses, domains, ranges,
%etc.)
%53C35:Global differential geometry- Symmetric spaces
\end{abstract}

%\newpage

\section{Introduction}

In Lie algebra, the classification theory of semisimple Lie algebras are well-established by Cartan, Killing and some others.   On the other hand, the classification theory of solvable or nilpotent Lie algebras are much less developed.     For nilpotent Lie algebras, there are only very few non-isomorphic nilpotent Lie algebras of dimension $n \le 6$.   Seeley \cite{Se93} and Gong \cite{Go98} classified nilpotent Lie algebras of dimension $7$.   According to \cite{Go98}, there are $119$ classes of indecomposable nilpotent Lie algebras of dimension $7$ over complex numbers, and $24$ additional classes over real numbers.   There are some partial results for dimension $8$.  The task for classifying all nilpotent Lie algebras  becomes unrealistic when the dimension $n$ reaches $9$.   For instance, according to \cite{TKK00}, there are $24,168$ non-isomorphic $9$-dimensional Lie algebras with a maximal Abelian ideal of dimension $7$.   Thus, for dimension $n \ge 8$, people turn their efforts to classify some special classes of nilpotent Lie algebras.

An important class of nilpotent Lie algebras is the class of two-step nilpotent Lie algebras.   The structure of Heisenberg algebras, a special class of two-step nilpotent Lie algebras, is well-known and simple.   As Heisenberg algebras are two-step nilpotent Lie algebras with centers of dimension $1$, people considered two-step nilpotent Lie algebras with centers of dimension $p \ge 2$.

In 1973, Gauger \cite{Ga73} studied two-step nilpotent Lie algebras.   He called them metabelian Lie algebras, and proposed a generator-relation method and a duality theory for two-step nilpotent Lie algebras.  By these, he classified two-step nilpotent Lie algebras for $n=8$, and proved that for each $n \ge 9$, there are infinitely many non-isomorphic two-step nilpotent Lie algebras.   In 1999, Galitski and Timashev \cite{GT99} continued this approach and classified two-step nilpotent Lie algebras for $n=9$.    %Their results have not excluded decomposable cases and given structure constants explicitly.   Hence, their classifications are somewhat incomplete.

In 2011, Ren and Zhu \cite{RZ11} adopted a different approach to attack this problem.
In 1971, Leger and Luks \cite{LL71} showed that the related sequences of a maximal torus minimal system of generators of a two-step nilpotent Lie algebra is an invariant of that algebra if the dimensions of the root spaces of the maximal torus are all $1$.   Using this result, Ren and Zhu \cite{RZ11} classified $8$-dimensional two-step nilpotent Lie algebras with a $2$-dimensional center.  They showed that there are at most five non-isomorphic indecomposable two-step nilpotent Lie algebras in this case.
On the other hand, Theorem 7.14 of \cite{Ga73} indicated that in the same case there are eleven non-isomorphic two-step nilpotent Lie algebras.   The difference is:  Ren and Zhu \cite{RZ11} excluded six decomposable two-step nilpotent Lie algebras.   See Example 1 in Section 3 of this paper.  It is also an important issue to identify a given Lie algebra being decomposable or indecomposable \cite{SW14}.  Hence, in a certain sense, Ren and Zhu \cite{RZ11} improved the result of Gauger \cite{Ga73} in this case.

Yan and Deng \cite{YD13}, Xia and Ren \cite{XR10}, Wang and Ren \cite{WR11}, Ren and Zhu \cite{RZ17, RZ17a} continued the approach of Ren and Zhu \cite{RZ11}, and classified $8$-dimensional and some $9$-dimensional two-step nilpotent Lie algebras.   Theorem 7.22 of \cite{Ga73} estimated that there are at most $42$ $8$-dimensional non-isomorphic two-step nilpotent Lie algebras with a $3$-dimensional center.   Galitski and Timashev \cite{GT99} indicated there are twelve non-isomorphic two-step nilpotent Lie algebras in this case, while Yan and Deng \cite{YD13} showed that there are at most eleven non-isomorphic indecomposable two-step nilpotent Lie algebras in this case.   Yan and Deng \cite{YD13} further excluded one decomposable algebra in this case.    Thus, Yan and Deng \cite{YD13} also improved the results of Gauger \cite{Ga73}, Galitski and Timashev \cite{GT99}, in this case.   See Example 4 in Section 3 of this paper.

%In 2017, Ren and Zhu \cite{RZ17a} further improved the result of Yan and Deng \cite{YD13}.  They pointed out that to show these eleven two-step nilpotent Lie algebras are non-isomorphic, the rank of a nilpotent Lie algebra \cite{GK96} may be used as the first invariant to distinguish them.   Then some two-step nilpotent Lie algebras for whom the dimensions of the root spaces of the maximal torus are not all $1$ may be distinguished.   Finally, the lemma of Leger and Luks \cite{LL71} may be used.

In the next section, we briefly describe Gauger's approach by the notion of free two-step nilpotent Lie algebras \cite{CFS18, GKT02}.

In Section 3, we define generating hypergraph for a two-step nilpotent Lie algebra.  The two-step nilpotent Lie algebra is decomposable if and only if its generating hypergraph is not connected under certain bases.  We use this property to identify some two-step nilpotent Lie algebras being decomposable or indecomposable.

In \cite{RZ11, RZ17, RZ17a, WR11, XR10, YD13}, a  two-step nilpotent Lie algebra of dimension $8$ or $9$ is denoted as $N^{n, p}_i$, where $n$ is the dimension of the two-step nilpotent Lie algebra, $p$ is the dimension of its center, $i$ is an auxiliary index without special meanings.  This causes some confusion and is not convenient for readers.    For instance, both \cite{YD13} and \cite{RZ17a} classified two-step nilpotent Lie algebra of $n=8$ and $p=3$ with different methods, and claimed that there are eleven non-isomorphic two-step nilpotent Lie algebras of $n=8$ and $p=3$, and denoted them as $N^{8, 3}_i$ for $i = 1, \cdots 11$.   However, the two algebras denoted as $N^{8, 3}_i$ for the same $i$ in \cite{YD13} and \cite{RZ17a} are not the same.  As the bases given in \cite{YD13} and \cite{RZ17a} are also different, it is not easy for readers to identify which two algebras in \cite{YD13} and \cite{RZ17a} are the same.
Hence, in Section 4, we introduce a conventional nomenclature for two-step nilpotent Lie algebras of dimension $n = 8, 9$, classified in \cite{RZ11, RZ17, RZ17a, WR11, XR10, YD13}.   We denote such two-step nilpotent Lie algebras as $T^{n, p}_r$, where $r$ is the rank of the two-step nilpotent Lie algebra, if there is only one non-isomorphic two-step nilpotent Lie algebra for such a value of $n, p$ and $r$.   Otherwise, we denote them as $T^{n, p}_{r, i}$.   In this way, the third index $r$ has a concrete meaning as the rank is the most important variant of a two-step nilpotent Lie algebra after $n$ and $p$.   The fourth index $i$ takes a larger value if the dimensions of some root spaces are greater than one.  When the root spaces are all of dimension one, the value of $i$ follows the lexicographic order of $H$-msg related sequences.

In doing this, we found that in several cases of \cite{RZ11, RZ17, RZ17a, YD13}, where there are no detailed proofs of the non-isomorphic property of those two-step nilpotent Lie algebras, when their dimensions, the dimensions of their centers, their ranks, and their $H$-msg related sequences are all the same.
In Section 5, beside the generating hypergraph, we define generator graph for a two-step nilpotent Lie algebra,  and introduce several new invariants associated with these two graphs.  By using these new invariants, we determine the values of the fourth index $i$ when the root spaces are all of dimension one, and the lexicographic order of $H$-msg related sequences is the same.

Ren and Zhu \cite{RZ11} showed that any indecomposable two-step nilpotent Lie algebra of $n=8$ and $p=2$ is isomorphic to one of the five two-step nilpotent Lie algebras.  In Section 6, by using the tool of generating hypergraphs, we give a direct proof that these five two-step nilpotent Lie algebras are all indecomposable.  Our method can be extended to the other cases.

Furthermore, by using the tool of generating hypergraphs, we identify some decomposable two-step nilpotent Lie algebras in dimension $9$ in Section 7.

The Lie algebras studied in this paper are over the complex number field.

\section{Free Two-Step Nilpotent Lie Algebras and Duality}

As in \cite{Ga73},  the notation $A \cong B$ is for isomorphic algebras $A$ and $B$.   We briefly describe Gauger's approach by the notion of free two-step nilpotent Lie algebras \cite{CFS18, GKT02}.

Suppose that $\{ \uu_1, \cdots, \uu_q \}$ is a set of generators, where $q \ge 2$.  Let $U$ be the space spanned by $\{ \uu_1, \cdots, \uu_q \}$.
Let $V = $ Span$\{ (\uu_i, \uu_j) : 1 \le i < j \le q \}$.     Then $V$ is a space of dimension $q(q-1)/2$.      Let $N^q = U \oplus V$.
For any $\uu_i, \uu_j$ with $1 \le i, j \le n$, and  $\vv_r, \vv_s \in V$, define

1) $[\uu_i, \uu_j] = (\uu_i, \uu_j) = - [\uu_j, \uu_i]$;

2) $[\uu_i, \vv_r] = 0$ and $[\vv_r, \vv_s] = 0$.

Then $N^q$ is a free two-step nilpotent Lie algebra of dimension $q(q+1)/2$ \cite{CFS18, GKT02}.

Suppose that $I$ is a subspace of $V$, with dimension $p$.  By the above properties of $N^q$,  $I$ is an ideal of $N^q$.   Let $N = N^q /I$.   Then $N$ is a two-step nilpotent Lie algebra with $q$ generators $\{ \x_i = \uu_i + I : 1 \le i \le q \}$.   Denote $D(N) = [N, N]$.   Then $D(N) = C(N)$ is the center of $N$, with dimension $\bar p = q(q-1)/2 - p$.    The dimension of $N$ is $q(q+1)/2 - p$.   Gauger \cite{Ga73} called $I$ the relation set of $N$.

The following proposition is a simple form of Proposition 1.6 of \cite{Ga73}.

\begin{proposition} \label{p1}
Let $I, J \subset V$.   Then $N^q/I \cong N^q/J$ if and only if there is an automorphism $\theta$ of $N^q$ such that $\theta(I) = J$.
\end{proposition}

On the other hand, by Theorem 2.1 of \cite{Ga73}, every two-step nilpotent Lie algebra with $q$-generators is of the type $N^q/I$, where $I$ is a proper subspace of $V$.
In this way, the classification of two-step nilpotent Lie algebras of $q$-generators and a $\bar p$-dimensional center is equivalent to classification of all factor algebras $N^q/I$ with a $p$-dimensional relation set $I$.

We now describe the duality theory of Gauger.   Suppose that there is an inner product defined on $V$.  For a subspace $I$ of $V$, there is an orthogonal complement subspace $I^\perp$ of $I$ in $V$.    Let  $N^\perp = N^q/I^\perp$.   We call $N^\perp$ the dual of $N$.     By Theorem 3.2 of \cite{Ga73}, we have the following theorem.

\begin{theorem} \label{t1}

(i) $N \cong (N^\perp)^\perp$;

(ii) $N_1 \cong N_2$ if and only if $N_1^\perp \cong N_2^\perp$;

(iii) if dim $N = q + \bar p$, then dim $N^\perp = q + p$.
\end{theorem}

By this duality theory, the classification of two-step nilpotent Lie algebras of $q$-generators and a $p$-dimensional center is equivalent to classification of all factor algebras $N^q/I$ with a $p$-dimensional relation set $I$.

Hence, $q$ and $p$ are two most important invariants of $N^\perp$, while $n=q+p$ is the dimension of $N^\perp$.  In \cite{GT99}, $(q, p)$ is called the signature of $N^\perp$.

For each $n$, the value of $p$ has an upper limit.

\begin{proposition}   \label{p2}
Let $N$ be an indecomposable $n$-dimensional two-step nilpotent Lie algebra, and its center have dimension $p$.
Then
\begin{equation} \label{e1}
1 \le p \le n +{1 \over 2} - \sqrt{2n + {1 \over 4}}.
\end{equation}
\end{proposition}
{\bf Proof} From the above discussion on $q$ and $p$, we have
$$n \le q + {q(q-1) \over 2} = {q(q+1) \over 2}.$$
This results (\ref{e1}). \qed

\section{Generating Hypergraph and Decomposability of Two-Step Nilpotent Lie Algebras}

Suppose that we have a two-step nilpotent Lie algebra $N$ with $q$-generators $\hat X = \{ \x_1, \cdots, \x_q \}$, and a center $I$, with a basis $\hat I = \{ \y_1, \cdots, \y_p \}$.  Let
$X=$ Span $\hat X$.    Then $N = X \oplus I$.   Assume that $q, p \ge 2$.    The dimension of $N$ is $n = q+p$.

 We now define {\bf generating hypergraph} $G$ for $N$.   The generating hypergraph $G$ is a bipartite hypergraph.  For basic knowledge of hypergraphs, see \cite{Br13}.   The bipartite hypergraph $G$ has two vertex sets $\hat X$ and $\hat I$.   If there is a definition of a Lie bracket operation
$$[\x_i, \x_j] = \sum_{k=1}^p \alpha_k \y_k,$$
where not all $\alpha_k$ are zero,
then $G$ has a multi-vertex hyper-edge $(\x_i, \x_j; \y_k : \alpha_k \not = 0 )$.   For each $k = 1, \cdots, q$, we assume that there is at least one three-vertex hyper-edge $(\x_i, \x_j, \y_k)$.  Otherwise, we may always make a linear transformation of the basis of $I$ to reach this assumption.

%The following proposition holds clearly.

\begin{proposition} \label{p3}
The two-step nilpotent Lie algebra $N$ is decomposable if and only if there are bases $\hat X$ and $\hat I$ such that the generating hypergraph $G$ is not connected.
\end{proposition}
{\bf Proof} Suppose that there are bases $\hat X$ and $\hat I$ such that the generating hypergraph $G$ is not connected.   Then the vertex set of $G$ can be partitioned to two nonempty  subsets $\{ \x_1, \cdots, \x_k, \y_1, \cdots, \y_l \}$ and  $\{ \x_{k+1}, \cdots, \x_q, \y_{l+1}, \cdots, \y_p \}$ such that these two vertex sets are not connected in $G$.  We see that $N$ is the direct sum of two two-step nilpotent Lie algebras $N_1$ and $N_2$, where $\{ \x_1, \cdots, \x_k \}$ is the generator set of $N_1$,  $\{ \y_1, \cdots, \y_l \}$ is a base of the center of $N_1$, $\{ \x_{k+1}, \cdots, \x_q \}$ is the generator set of $N_2$,  $\{ \y_{l+1}, \cdots, \y_p \}$ is a base of the center of $N_2$.

On the other hand, suppose that $N$ is the direct sum of two two-step nilpotent Lie algebras $N_1$ and $N_2$.  Take the generator sets, and bases of the centers of $N_1$ and $N_2$, we may have bases $\hat X$ and $\hat I$ such that the generating hypergraph $G$ formed from these bases is not connected.  \qed

However, this property is basis-dependent.   Let $N$ be defined by $[\x_1, \x_2] = \y_1$ and $[\x_3, \x_4] = \y_2$.  Then $N$ is decomposable by this proposition.   Now let $\bar \x_2 = \x_2 + \x_3$, $\bar \x_3 = \x_2 - \x_3$, and use $\bar \x_2$ and $\bar \x_3$ to replace $\x_2$ and $\x_3$.  Then $N$ is defined by $[\x_1, \bar \x_2] = \y_1$, $[\x_1, \bar \x_3] = \y_1$,
$[\bar \x_2, \x_4] = \y_2$ and $[\x_4, \bar \x_3] = \y_2$.   The resulted generating hypergraph is connected.

This proposition provides us a tool to exclude decomposable two-step nilpotent Lie algebras from the classification.  For more knowledge about connectivity of a hypergraph, the related quantities such as algebraic connectivity and analytical connectivity, and their computational methods, see \cite{QL17}.

{\bf Example 1} Theorem 7.14 of \cite{Ga73} claimed that every $6$-generator and $2$-relation two-step nilpotent Lie algebra $N$ is isomorphic to exactly one of $N^6/I_j$ for $j = 1, \cdots, 11$:
$$I_1 = {\rm Span}\{ [\uu_1, \uu_2] + [\uu_5, \uu_6], [\uu_3, \uu_4] + [\uu_5, \uu_6] \},$$
$$I_2 = {\rm Span}\{ [\uu_1, \uu_2] + [\uu_3, \uu_4], [\uu_5, \uu_6] \},$$
$$I_3 = {\rm Span}\{ [\uu_1, \uu_4] + [\uu_2, \uu_3], [\uu_2, \uu_4] + [\uu_5, \uu_6] \},$$
$$I_4 = {\rm Span}\{ [\uu_1, \uu_4] + [\uu_2, \uu_3] + [\uu_5, \uu_6], [\uu_2, \uu_4]  \},$$
$$I_5 = {\rm Span}\{ [\uu_1, \uu_6] + [\uu_2, \uu_5] + [\uu_3, \uu_4], [\uu_2, \uu_6] + [\uu_3, \uu_5] \},$$
$$I_6 = {\rm Span}\{ [\uu_1, \uu_2], [\uu_3, \uu_4] \},$$
$$I_7 = {\rm Span}\{ [\uu_1, \uu_4] + [\uu_2, \uu_3], [\uu_2, \uu_4] \},$$
$$I_8 = {\rm Span}\{ [\uu_1, \uu_2] + [\uu_5, \uu_6], [\uu_4, \uu_6] \},$$
$$I_9 = {\rm Span}\{ [\uu_5, \uu_6], [\uu_4, \uu_6] \},$$
$$I_{10} = {\rm Span}\{ [\uu_1, \uu_3] + [\uu_4, \uu_6], [\uu_2, \uu_3] + [\uu_5, \uu_6] \},$$
$$I_{11} = {\rm Span}\{ [\uu_2, \uu_6] + [\uu_3, \uu_5], [\uu_3, \uu_6] + [\uu_4, \uu_5] \}.$$
Here, $U = \{ \uu_1, \cdots, \uu_6 \}$ for $N^6$.  Now, we check $I_6$.  Then
$$I_6^\perp = {\rm Span}\{ [\uu_i, \uu_j] : 1 \le i < j \le 6, (1, 2) \not = (i, j) \not = (3, 4) \},$$
and
$$N^6/I_6^\perp = {\rm Span} \{ \x_1, \cdots, \x_6, \y_1 = [\x_1, \x_2], \y_2 = [\x_3, \x_4] \}.$$
By Proposition \ref{p3}, we see that $N^6/I_6^\perp$ is decomposable.   Similarly, we may conclude that $N^6/I_j^\perp$ are decomposable for $j = 2, 7, 8, 9, 11$.
By the discussion in the last section, we see that $N^6/I_{10}^\perp = N^{8, 2}_2$ and $N^6/I_j^\perp = N^{8, 2}_j$ for $j = 1, 3, 4, 5$, where $N^{8, 2}_j$ for $j = 1, \cdots, 5$, are five
 non-isomorphic two-step nilpotent Lie algebras given by Ren and Zhu in \cite{RZ11}.  In the next section, we will express $N^{8, 2}_j$ for $j = 1, \cdots, 5$, explicitly.   We will make a conventional nomenclature for them in the next two sections.

 {\bf Example 2}  Theorem 7.12 of \cite{Ga73} claimed that every $4$-generator and $2$-relation two-step nilpotent Lie algebra $N$ is isomorphic to exactly one of $N^4/I_j$ for $j = 1, 2, 3$:
$$I_1 = {\rm Span}\{ [\uu_1, \uu_2], [\uu_3, \uu_4] \},$$
$$I_2 = {\rm Span}\{ [\uu_1, \uu_4] + [\uu_2, \uu_3], [\uu_2, \uu_4] \},$$
$$I_3 = {\rm Span}\{ [\uu_2, \uu_4], [\uu_3, \uu_4] \}.$$
Here, $U = \{ \uu_1, \cdots, \uu_4 \}$ for $N^4$.  By the discussion in the last section, we see that $N^4/I_1 = N^{8, 4}_1$ and $N^4/I_2 = N^{8, 4}_3$ and $N^4/I_3 = N^{8, 4}_2$,  where $N^{8, 4}_j$ for $j = 1, 2, 3$, are three
 non-isomorphic two-step nilpotent Lie algebras given by Yan and Deng for $n=8, p=4$ in \cite{YD13}.  Also see \cite{XR10}.  In the next section, we will express $N^{8, 4}_j$ for $j = 1, 2, 3$, explicitly,  and  will make a conventional nomenclature for them.

 {\bf Example 3}  Theorem 5.2 of \cite{Ga73} claimed that every $4$-generator and $1$-relation two-step nilpotent Lie algebra $N$ is isomorphic to exactly one of $N^4/I_j$ for $j = 1, 2$:
$$I_1 = {\rm Span}\{ [\uu_1, \uu_2] \},$$
$$I_2 = {\rm Span}\{ [\uu_1, \uu_2] + [\uu_3, \uu_4] \}.$$
Here, $U = \{ \uu_1, \cdots, \uu_4 \}$ for $N^4$.  By the discussion in the last section, we see that $N^4/I_1 = N^{9, 5}_1$ and $N^4/I_2 = N^{9, 5}_2$,  where $N^{9, 5}_j$ for $j = 1, 2$, are two
 non-isomorphic two-step nilpotent Lie algebras given by Wang and Ren for $n=9, p=5$ in \cite{WR11}.   In the next section, we will express $N^{9, 5}_j$ for $j = 1, 2$, explicitly,  and  will make a conventional nomenclature for them.

 {\bf Example 4} Theorem 7.22 of \cite{Ga73} said that there are at most $42$ non-isomorphic non-isomorphic $5$-generator, $3$-relation two-step nilpotent Lie algebras. This corresponds to the case that $n=8$ and $p=3$. In Tables 2 and 6 of \cite{GT99}, Galitski and Timashev indicated there are only twelve non-isomorphic two-step nilpotent Lie algebras in this case.    The No. 91 algebra of Table 2 of \cite{GT99} is defined by
 $$[\x_1, \x_4] = \y_2, [\x_1, \x_5] = \y_3, [\x_2, \x_3] = \y_1.$$
 Consider the generating hypergraph $G$ of this two-step nilpotent Lie algebra.  Then the vertex set of $G$ can be partitioned to two sets $\{ \x_1, \x_4, \x_5, \y_2, \y_3 \}$ and $\{ \x_2, \x_3, \y_1 \}$.  These two vertex sets are not connected in $G$.   Thus, by Proposition \ref{p3}, this two-step nilpotent Lie algebra is decomposable.  Otherwise, the No. 87, 90, 78, 80, 92, 88, 75, 72, 84, 68, 62 algebras of Table 2 of \cite{GT99} correspond $N^{8,3}_j$ in \cite{RZ17a} for $j=1, \cdots, 11$, respectively.  In the next section, we will express $N^{8, 3}_j$ for $j = 1, \cdots, 11$, explicitly.   We will make a conventional nomenclature for them in the next two sections.   Note that the eleven two-step nilpotent Lie algebras in \cite{RZ17a} are the same as the eleven two-step nilpotent Lie algebras in \cite{YD13}, with different orders and expressions.

\section{Rank, Dimension of Root Space, and $H$-msg Related Sequence}

Assume that $ p \ge 2$.  By (\ref{e1}), for $n = 8$, $p=2, 3$ and $4$; for $n = 9$, $p=2, 3, 4$ and $5$.    Ren and Zhu \cite{RZ11, RZ17a}, Xia and Ren \cite{XR10}, and Yan and Deng \cite{YD13} made the classification of two-step nilpotent Lie algebras of dimension $8$.   Wang and Ren \cite{WR11}, and Ren and Zhu \cite{RZ17} made the classification of two-step nilpotent Lie algebras of dimension $9$ with center dimension $5$ and $2$.

Let $N$ be a two-step nilpotent Lie algebra, $I$ be its center, $n$, $p$ and $q$ be as defined above.  Assume that $n \ge 8$ and $p \ge 2$. Suppose that $\{ \x_1, \cdots, \x_q \}$ is a minimal system of generators of $N$.   By \cite{RZ11, RZ17, YD13}, the related set of $\x_i$ is defined to be the set $G(\x_i) = \{ \x_j : [\x_i, \x_j] \not = \0\}$, the number $p_i = |G(\x_i)|$ is called related number of $\x_i$, the $q$-tuple of integers $(p_1, \cdots, p_q)$ is called the related sequence of $\{ \x_1, \cdots, \x_q \}$ \cite{Fa73}.   We will order $\x_1, \cdots, \x_q$ such that $p_i \le p_{i+1}$.

 Denote the set of all derivations of $N$ by Der$(N)$.  A maximal torus $H$ of $N$ is a maximal Abelian subalgebra of Der$(N)$, which consists of semi-simple linear transformation.  Then $N$ can be decomposed into a direct sum of root spaces for $H$:
 $N = \sum_{\beta \in H^*} N_\beta$, where $H^*$ is the dual space of $H$, and
 $$N_\beta = \{ \x \in N : h(\x) = \beta(h) \x, \forall h \in H \}.$$
 If $H$ is a maximal torus on $N$, let $\Delta = \{ \beta \in H^* : N_\beta \not = 0 \}$ be the root system of $N$ associated with $H$.

 For a two-step nilpotent Lie algebra $N$, a maximal torus always exists \cite{RZ11}.
 By Mostow's theorem \cite[Theorem 4.1]{Mo56}, all maximal tori of a nilpotent Lie algebra have the same dimension $r$.   Hence, this dimension $r$ is an invariant of the nilpotent Lie algebra $N$, called the rank of $N$, and denoted as $r=$ Rank$(N)$ \cite{GK96}.

 Let $H$ be a maximum torus of $N$.  A minimal system of generators consisting of root vectors for $H$ is called an $H$-msg of $N$ \cite{Sa82}.

 A minimal system of generators is called a $(p_1, \cdots, p_q)$-msg if its related sequence is $(p_1, \cdots, p_q)$.  It is called a $(p_1, \cdots, p_q)$-$H$-msg if it is also an $H$-msg \cite{RZ11}.

 A key lemma used in \cite{RZ11, RZ17, RZ17a, WR11, XR10, YD13} was originally from \cite{LL71}, see \cite{YD13}.   By \cite{LL71}, a nilpotent Lie algebra $N$ is called quasi-cyclic if $N$ has a subspace $U$ such that $N = U \oplus U^1 \oplus \cdots, \oplus U^k$, where $U^0 = U, U^i = [U, U^{i-1}]$.    Then a two-step Nilpotent Lie algebra, is quasi-cyclic.

 \begin{lemma} \label{l1} \cite{LL71}
 Let $N$ be a quasi-cyclic nilpotent Lie algebra, and $\{ \x_1, \cdots, \x_q \}$ an $H_1$-msg of $N$, $\{ \z_1, \cdots, \z_q \}$ an $H_2$-msg of $N$.
 Then there exists an automorphism $\theta$ of $N$ such that
 $$\left(\z_1,\cdots,\z_q\right)^\top = A\left(\theta(\x_1),\cdots,\theta(\y_q)\right)^\top,$$
 where  $\left(\z_1,\cdots,\z_q\right)^\top$ is the transpose of the matrix $\left(\z_1,\cdots,\z_q\right)$, and $A$ is a $q \times q$ matrix.   In particular, if for any $i$, the dimension of the root space associated with $\x_i$ is $1$, then $A$ is a monomial matrix (i.e., each row or each column of $A$ has exactly one nonzero entry).
 \end{lemma}

 In general, the related sequence of an $H$-msg is not necessarily an invariant \cite{RZ11}.
 For a two-step nilpotent Lie algebra $N$, the related sequence $(p_1, \cdots, p_q)$ of an $H$-msg is an invariant if the dimensions of all the root spaces of $H$ are of dimension $1$ by Lemma \ref{l1}.

We denote such two-step nilpotent Lie algebras as $T^{n, p}_r$, where $r$ is the rank of the two-step nilpotent Lie algebra, if there is only one non-isomorphic two-step nilpotent Lie algebra for such a value of $n, p$ and $r$.   Otherwise, we denote them as $T^{n, p}_{r, i}$.   The fourth index $i$ takes a larger value if the dimensions of some root spaces are greater than one.  When the root spaces are all of dimension one, the value of $i$ follows the lexicographic order of $H$-msg related sequences.  When the root spaces are all of dimension one, and the lexicographic order of $H$-msg related sequences is the same, the value of $i$ will be determined in the next section.

\bigskip

We now consider the case that $n=8$.

By \cite{RZ11, YD13}, for $n=8$ and $p=2$, there are five possibly non-isomorphic two-step nilpotent Lie algebras $N^{8, 2}_i$ for $i = 1, \cdots, 5$.

The two-step nilpotent Lie algebra $N^{8, 2}_1$ is defined by
$$[\x_1, \x_2] = \y_1, [\x_3, \x_4] = \y_2, [\x_5, \x_6] = \y_1 + \y_2.$$
It has an $H$-msg related sequence $(1, 1, 1, 1, 1, 1)$.  Its rank $r$ is $4$.

The two-step nilpotent Lie algebra $N^{8, 2}_2$ is defined by
$$[\x_5, \x_2] = [\x_6, \x_1] = \y_1, [\x_5, \x_3] = [\x_6, \x_4] = \y_2.$$
The two-step nilpotent Lie algebra $N^{8, 2}_3$ is defined by
$$[\x_1, \x_2] = [\x_6, \x_5] = \y_1, [\x_3, \x_6] = [\x_5, \x_4] = \y_2.$$
The two-step nilpotent Lie algebra $N^{8, 2}_4$ is defined by
$$[\x_1, \x_2] = [\x_3, \x_6] = [\x_5, \x_4] = \y_1, [\x_6, \x_5] = \y_2.$$
Then the  two-step nilpotent Lie algebras $N^{8, 2}_i$ for $i = 2, 3, 4$ have an $H$-msg related sequence $(1, 1, 1, 1, 2, 2)$.   They also have $r = 4$.

The two-step nilpotent Lie algebra $N^{8, 2}_5$ is defined by
$$[\x_1, \x_6] = [\x_3, \x_4] = [\x_5, \x_2] = \y_1, [\x_6, \x_3] = [\x_4, \x_5] = \y_2.$$
It has an $H$-msg related sequence $(1, 1, 2, 2, 2, 2)$.  Its rank $r$ is $3$.   This shows that $N^{8, 2}_5$ is not isomorphic from $N^{8, 2}_i$ with $i = 1, 2, 3, 4$.  Then we denote it as $T^{8, 2}_3$.

The dimensions of the root spaces of $N^{8, 2}_i$ with $i = 1, 2, 3, 4$ are of $1$.   Under this condition, their $H$-msg related sequences are invariants.  Since $(1, 1, 1, 1, 1, 1) \prec (1, 1, 1, 1, 2, 2)$, $N^{8, 2}_1$ is not isomorphic from $N^{8, 2}_i$ with $i = 2, 3, 4$, and we denote it as $T^{8, 2}_{4, 1}$,    We tentative;y denote $N^{8, 2}_i$ with $i = 2, 3, 4$ as $T^{8, 2}_{4, i}$ with $i = 2, 3, 4$, and do not specify the exact value of $i$ for each of them at this moment.
In the next section, we will show that $N^{8, 2}_i$ with $i = 2, 3, 4$ are mutually non-isomorphic and determine the exact value of $i$ in $T^{8, 2}_{4, i}$ with $i = 2, 3, 4$ for these three two-step nilpotent Lie algebras.

By \cite{RZ17a}, for $n=8$ and $p=3$, there are eleven possibly non-isomorphic two-step nilpotent Lie algebras $N^{8, 3}_i$ for $i = 1, \cdots, 11$.

The two-step nilpotent Lie algebra $N^{8, 3}_1$ is defined by
$$[\x_1, \x_2] = [\x_3, \x_4] = \y_1, [\x_3, \x_5] = \y_2, [\x_4, \x_5] = \y_3.$$
The two-step nilpotent Lie algebra $N^{8, 3}_2$ is defined by
$$[\x_1, \x_5] = [\x_4, \x_2] = \y_1, [\x_5, \x_3] = \y_2, [\x_3, \x_4] = \y_3.$$
The two-step nilpotent Lie algebra $N^{8, 3}_3$ is defined by
$$[\x_5, \x_3] = [\x_3, \x_4] = \y_1, [\x_1, \x_5] = \y_2, [\x_2, \x_4] = \y_3.$$
The two-step nilpotent Lie algebra $N^{8, 3}_4$ is defined by
$$[\x_1, \x_5] = [\x_3, \x_4] = \y_1, [\x_5, \x_3] = \y_2, [\x_2, \x_4] = \y_3.$$
Then the two-step nilpotent Lie algebras $N^{8, 3}_i$ for $i = 1, 2, 3, 4$ have an $H$-msg related sequence $(1, 1, 2, 2, 2)$.   By \cite{RZ17a}, they have $r = 4$.

The two-step nilpotent Lie algebra $N^{8, 3}_5$ is defined by
$$[\x_1, \x_4] = [\x_5, \x_2] = \y_1, [\x_3, \x_5] = \y_2, [\x_4, \x_5] = \y_3.$$
Then it has an $H$-msg related sequence $(1, 1, 1, 2, 3)$.  %It also has $s = 8$.
By \cite{RZ17a},  it has $r = 4$.

The two-step nilpotent Lie algebra $N^{8, 3}_8$ is defined by
$$[\x_1, \x_5] = [\x_3, \x_4] = \y_1, [\x_3, \x_5] = [\x_2, \x_4] = \y_2, [\x_4, \x_5] = \y_3.$$
Thus it has an $H$-msg related sequence $(1, 1, 2, 3, 3)$.   By \cite{RZ17a}, it has $r = 3$.

The two-step nilpotent Lie algebra $N^{8, 3}_6$ is defined by
$$[\x_1, \x_2] = [\x_5, \x_4] = \y_1, [\x_2, \x_5] = [\x_4, \x_3] = \y_2, [\x_3, \x_5] = \y_3.$$
The two-step nilpotent Lie algebra $N^{8, 3}_7$ is defined by
$$[\x_1, \x_2] = [\x_5, \x_3] = \y_1, [\x_2, \x_5] = [\x_5, \x_4] = \y_2, [\x_3, \x_4] = \y_3.$$
The two-step nilpotent Lie algebra $N^{8, 3}_9$ is defined by
$$[\x_1, \x_5] = [\x_3, \x_2] = \y_1, [\x_3, \x_5] = [\x_2, \x_4] = \y_2, [\x_4, \x_5] = \y_3.$$
Then the two-step nilpotent Lie algebras $N^{8, 3}_i$ for $i=6, 7, 9$ have an $H$-msg related sequence $(1, 2, 2, 2, 3)$.    By \cite{RZ17a}, they also have  $r = 3$.

The two-step nilpotent Lie algebra $N^{8, 3}_{10}$ is defined by
$$[\x_1, \x_2] = [\x_3, \x_5] = \y_1, [\x_2, \x_3] = [\x_4, \x_5] = \y_2, [\x_1, \x_5] = \y_3.$$
Then it has an $H$-msg related sequence $(2, 2, 2, 2, 2)$.   By \cite{RZ17a}, it also has $r = 3$.

The two-step nilpotent Lie algebra $N^{8, 3}_{11}$ is defined by
$$[\x_1, \x_5] = [\x_4, \x_2] = \y_1, [\x_1, \x_4] = [\x_5, \x_3] = \y_2, [\x_4, \x_5] = [\x_2, \x_3] = \y_3.$$
Thus it has an $H$-msg related sequence $(2, 2, 2, 3, 3)$.    By \cite{RZ17a}, it has $r=2$.   Since it is the only one among these eleven two-step nilpotent Lie algebras with $r=2$, it is not isomorphic from $N^{8, 3}_i$ with $i = 1, \cdots, 10$.  Then, we denote it as $T^{8, 3}_2$.

In \cite{RZ17a}, in the proof of Theorem 1, it was said that the dimension of the maximal torus of $N^{8, 3}_7$, i.e., the rank of $N^{8, 3}_7$, is $4$.  This should be a typo.  In fact, in the latter part of the proof of Theorem 1 of \cite{RZ17a}, it was said that the dimensions of the maximal toruses of $N^{8, 3}_i$ with $i = 6, \cdots, 10$ are the same.  Dr. Zaili Yan, the first author of \cite{YD13}, also showed us that $N^{8, 3}_7$ of \cite{RZ17a} is isomorphic with $N^{8, 3}_{11}$ of \cite{YD13}, while the rank of $N^{8, 3}_{11}$ of \cite{YD13} is $3$.   See Section 6 for details.

 The ranks of $N^{8, 3}_i$ for $i = 1, \cdots, 5$ are all $4$.  The ranks of $N^{8, 3}_j$ for $j = 6, \cdots, 10$ are all $3$.   This shows that $N^{8, 3}_i$ and $N^{8, 3}_j$ are not isomorphic for any $i = 1, \cdots, 5$ and any $j = 6, \cdots, 10$.

The ranks of $N^{8, 3}_i$ for $i = 1, \cdots, 5$ are all $4$.   However, the root space of $N^{8, 3}_3$, related with $\x_3$, has dimension $2$, while all the root spaces of $N^{8, 3}_i$
with $i = 1, 2, 4, 5$, have dimension $1$.  This implies that $N^{8, 3}_3$ is not isomorphic from $N^{8, 3}_i$
with $i = 1, 2, 4, 5$, and its structure is more complicated.  Thus, we denote $N^{8, 3}_3$ as $T^{8, 3}_{4, 5}$.

Since all the root spaces of $N^{8, 3}_i$ with $i = 1, 2, 4, 5$, have dimension $1$, their $H$-msg related sequences are invariants.
The $H$-msg related sequence of $N^{8, 3}_5$ is $(1, 1, 1, 2, 3)$; while the $H$-msg related sequence of $N^{8, 3}_i$
with $i = 1, 2, 4$ is $(1, 1, 2, 2, 2)$.  This shows that $N^{8, 3}_5$ is not isomorphic from $N^{8, 3}_i$
with $i = 1, 2, 4$.   Since $(1, 1, 1, 2, 3) \prec (1, 1, 2, 2, 2)$, we denote $N^{8, 3}_5$ as $T^{8, 3}_{4, 1}$.  We will show that $N^{8, 3}_i$
with $i = 1, 2, 4$ are mutually non-isomorphic, and determine the value of $i$ in $T^{8, 3}_{4, i}$ with $i = 2, 3, 4$ for each of them, in the next section.

The ranks of $N^{8, 3}_j$ for $j = 6, \cdots, 10$ are all $3$.  However, the root space of $N^{8, 3}_7$, related with $\x_5$, has dimension $2$, while all the root spaces of $N^{8, 3}_i$
with $i = 6, 8, 9, 10$, have dimension $1$.  This implies that $N^{8, 3}_7$ is not isomorphic from $N^{8, 3}_i$
with $i = 6, 8, 9, 10$, and its structure is more complicated.  Thus, we denote $N^{8, 3}_7$ as $T^{8, 3}_{3, 5}$.

Since all the root spaces of $N^{8, 3}_i$ with $i = 6, 8, 9, 10$, have dimension $1$, their $H$-msg related sequences are invariants.  The $H$-msg related sequence of $N^{8, 3}_8$ is $(1, 1, 2, 3, 3)$.   The $H$-msg related sequences of $N^{8, 3}_i$ with $i = 6, 9$ are $(1, 2, 2, 2, 3)$.
The $H$-msg related sequence of $N^{8, 3}_{10}$ is $(2, 2, 2, 2, 2)$.  They are different.   This shows that $N^{8, 3}_8$, $N^{8, 3}_{10}$ and any of
$N^{8, 3}_i$ with $i = 6, 9$ are mutually non-isomorphic.  Since $(1, 1, 2, 3, 3) \prec (1, 2, 2, 2, 3) \prec (2, 2, 2, 2, 2)$, we denote $N^{8, 3}_8$
as $T^{8, 3}_{3, 1}$, $N^{8, 3}_{10}$ as $T^{8, 3}_{3, 4}$, respectively.    We will show that $N^{8, 3}_i$
with $i = 6, 9$ are mutually non-isomorphic, and determine the value of $i$ in $T^{8, 3}_{3, i}$ with $i = 2, 3$ for each of them, in the next section.

By \cite{YD13}, for $n=8$ and $p=4$, there are three possible non-isomorphic two-step nilpotent Lie algebras $N^{8, 4}_i$ for $i = 1, 2, 3$.

The two-step nilpotent Lie algebra $N^{8, 4}_1$ is defined by
$$[\x_1, \x_2] = \y_1, [\x_2, \x_3] = \y_2, [\x_3, \x_4] = \y_3, [\x_4, \x_1] = \y_5.$$
Thus it has an $H$-msg related sequence $(2, 2, 2, 2)$.    By \cite{YD13}, it has $r = 4$.

The two-step nilpotent Lie algebra $N^{8, 4}_2$ is defined by
$$[\x_2, \x_4] = \y_1, [\x_3, \x_4] = \y_2, [\x_2, \x_3] = \y_3, [\x_1, \x_4] = \y_5.$$
Thus it has an $H$-msg related sequence $(1, 2, 2, 3)$.    By \cite{YD13}, it also has $r = 4$.

The two-step nilpotent Lie algebra $N^{8, 4}_3$ is defined by
$$[\x_3, \x_4] = \y_1, [\x_1, \x_3] = \y_2, [\x_2, \x_4] = y_2, [\x_1, \x_4] = \y_3, [\x_2, \x_3] = \y_4.$$
Thus it has an $H$-msg related sequence $(2, 2, 3, 3)$.   By \cite{YD13}, it has $r = 3$.    This shows that it is not isomorphic from $N^{8, 4}_1$
and $N^{8, 4}_2$.   We thus denote it as $T^{8, 4}_3$.

The ranks of $N^{8, 4}_1$ and $N^{8, 4}_2$ are the same.  Their root spaces are all of dimension $1$.  Thus their $H$-msg related sequences are invariants.  Since their $H$-msg related sequences are different, they are not isomorphic.  Since $(1, 2, 2, 3) \prec (2, 2, 2, 2)$, we denote
$N^{8, 4}_2$ as $T^{8, 4}_{4, 1}$, and $N^{8, 4}_1$ as $T^{8, 4}_{4, 2}$, respectively.

\medskip

We now consider the case that $n=9$.  In the literature, for $n=9$ and $p \ge 2$, the cases $p = 2$ and $5$ are known \cite{RZ17, WR11}.

By \cite{RZ17}, for $n=9$ and $p=2$, there are five possible non-isomorphic two-step nilpotent Lie algebras $N^{9, 2}_i$ for $i = 1, \cdots, 5$.

The two-step nilpotent Lie algebra $N^{9, 2}_1$ is defined by
$$[\x_1, \x_2] = [\x_4, \x_5] = [\x_6, \x_7] = \y_1, [\x_3, \x_7] = \y_2.$$
The two-step nilpotent Lie algebra $N^{9, 2}_2$ is defined
$$[\x_2, \x_7] = [\x_4, \x_5] = \y_1, [\x_1, \x_3] = [\x_6, \x_7] = \y_2.$$
Then they have an $H$-msg related sequence $(1, 1, 1, 1, 1, 1, 2)$.  Their rank $r = 5$.

The two-step nilpotent Lie algebra $N^{9, 2}_3$ is defined by
$$[\x_1, \x_2] = [\x_3, \x_7] = [\x_5, \x_6] = \y_1, [\x_4, \x_6] = [\x_5, \x_7] = \y_2.$$
The two-step nilpotent Lie algebra $N^{9, 2}_4$ is defined by
$$[\x_7, \x_2] = [\x_4, \x_5] = [\x_6, \x_1] = \y_1, [\x_7, \x_3] = [\x_5, \x_6] = \y_2.$$
Thus they have an $H$-msg related sequence $(1, 1, 1, 1, 2, 2, 2)$.  Their rank $r = 4$

The two-step nilpotent Lie algebra $N^{9, 2}_5$ is defined by
$$[\x_1, \x_7] = [\x_3, \x_4] = [\x_5, \x_6] = \y_1, [\x_7, \x_3] = [\x_4, \x_5] = [\x_6, \x_2] = \y_2.$$
Thus it has an $H$-msg related sequence $(2, 2, 2, 2, 2, 1, 1)$.   It has rank $r = 3$.

Because of the difference of their ranks, $N^{9, 2}_5$, any one of $N^{9, 2}_3$ and $N^{9, 2}_4$,  and any one of $N^{9, 2}_1$ and $N^{9, 2}_2$,
are mutually non-isomorphic.   In particular, $N^{9, 2}_5$ is not isomorphic from the other four two-step nilpotent Lie algebras.   Hence, we denote
$N^{9, 2}_5$ as $T^{9, 2}_3$.

However, the ranks and $H$-msg related sequences of $N^{9, 2}_1$ and $N^{9, 2}_2$ are the same, and the ranks and $H$-msg related sequences of $N^{9, 2}_3$ and $N^{9, 2}_4$ are the same.  Thus, we cannot further distinguish them at this moment.   In the next section, we will show that $N^{9, 2}_1$ and $N^{9, 2}_2$ are not isomorphic, $N^{9, 2}_3$ and $N^{9, 2}_4$ are not isomorphic, and determine the value of $i = 1, 2$ in $T^{9, 2}_{5, i}$ for
$N^{9, 2}_1$ and $N^{9, 2}_2$, and the value of $i = 1, 2$ in $T^{4, 2}_{5, i}$ for
$N^{9, 2}_3$ and $N^{9, 2}_4$.

By \cite{WR11}, for $n=9$ and $p=5$, there are two possible non-isomorphic two-step nilpotent Lie algebras $N^{9, 5}_1$ and $N^{9, 5}_2$.

The two-step nilpotent Lie algebra $N^{9, 5}_1$ is defined by
$$[\x_1, \x_3] = \y_1, [\x_1, \x_4] = \y_2, [\x_2, \x_3] = \y_3, [\x_2, \x_4] = \y_4, [\x_3, \x_4] = \y_5.$$
Then it has an $H$-msg related sequence $(2, 2, 3, 3)$.   We may check that its rank $r =4$.

The two-step nilpotent Lie algebra $N^{9, 5}_2$ is defined by
$$[\x_1, \x_2] = \y_1, [\x_1, \x_3] = [\x_2, \x_4] = \y_2, [\x_1, \x_4] = \y_3, [\x_2, \x_3] = \y_4, [\x_3, \x_4] = \y_5.$$
Then it has an $H$-msg related sequence $(3, 3, 3, 3)$.   We may check that its rank $r =4$ too.

The ranks of these two two-step nilpotent Lie algebras are all $4$.    We also see that the dimensions of their root spaces are all $1$.   Thus, their $H$-msg related sequences are invariants.   Since $(2, 2, 3, 3) \prec (3, 3, 3, 3)$, they are not isomorphic.    We denote  denoted $N^{9, 5}_1$ as $T^{9, 5}_{4, 1}$, and $N^{9, 5}_2$ as $T^{9, 5}_{4, 2}$, respectively.

\section{Generator Graph and Some New Invariants}

In the last section, there are five cases, where the non-isomorphic property is not established and the value of the fourth index $i$ is undetermined. In this section, by introducing more new invariants, we complete the non-isomorphic property proof, and determine the values of the fourth index $i$ in our nomenclature.

 Recall the generating hypergraph $G$ for a two-step nilpotent Lie algebra $N$, defined in Section 3.  We assume that the basis of $I$ is chosen such that the sum of the sizes of the hyper-edges of $G$ is minimized.   This assumption is important, otherwise center related sequence and weighted center related sequence introduced below may not be invariants.   In the following cases, the generating graphs involved are all $3$-uniform hypergraphs, i.e., each hyper-edge has exactly three vertices $\x_i, \x_j, \y_k$.   Then this assumption is satisfied.

From $G$, we further define an ordinary graph $\hat G$, called the {\bf generator graph} of $N$,  The vertex set of $\hat G$ is $\hat X$.  If the Lie bracket operation $[\x_i, \x_j]$ is defined, then $(\x_i, \x_j)$ is an edge of $\hat G$.   The generating graph $G$ and the generator graph $\hat G$ will be helpful to calculate the related sequence, the related index, and some new invariants we will introduce below.

We now study the case that $n=8$ and $p=2$.   In that case, the value of the fourth index $i$ in our nomenclature is undetermined for $T^{8, 2}_{4, i}$ with $i = 1, 2, 3$, where the $H$-msg related sequences are the same.  By our discussion in the last section, the candidates for $T^{8, 2}_{4, i}$ with $i = 1, 2, 3$ are $N^{8, 2}_i$ with $i = 2, 3, 4$ in \cite{RZ11}.

In \cite{RZ11}, the five two-step nilpotent Lie algebras $N^{8, 2}_i$ with $i = 1, 2, 3, 4, 5$ are claimed to be mutually non-isomorphic by citing a lemma (Lemma 5 of \cite{RZ11}) that if the dimensions of all the root spaces are $1$, then the Lie bracket relations are reserved except some index exchanges and scalings.    This implies that the $H$-msg related sequences are invariants under this condition.   However, for $N^{8, 2}_i$ with $i = 2, 3, 4$, the $H$-msg related sequences are the same.  How can we distinguish them quantitatively?  We now introduce two more new invariants.

First, we  put the dimensions of the components of the generator graph $\hat G$ in the nondecreasing order to make a sequence, and call it the {\bf generator relation sequence}.  Under the condition of Lemma 5 of \cite{RZ11}, this sequence is clearly an invariant.    Now, the generator relation sequences of $N^{8, 2}_i$ with $i = 2, 3, 4$ in \cite{RZ11} are $(3, 3)$, $(2, 4)$ and $(2, 4)$, respectively.  This shows that $N^{8, 2}_2$ are non-isomorphic from $N^{8, 2}_i$ with $i = 3, 4$.

Second, for each element of the center $I$, we count the number of its related Lie bracket operations, and put them in the nondecreasing order.  We call this sequence the {\bf center related sequence}.   By definition, we may use the generating graph $G$ to calculate this sequence.   Under the condition of Lemma \ref{l1} and the condition that $G$ is a $3$-uniform hypergraph, this sequence is also an invariant.   Now, the center related sequences of $N^{8, 2}_i$ with $i = 2, 3, 4$ in \cite{RZ11} are $(2, 2)$, $(2, 2)$ and $(1, 3)$, respectively.  This shows that $N^{8, 2}_4$ are non-isomorphic from $N^{8, 2}_i$ with $i = 2, 3$.

Then we have the following theorem.

\begin{theorem} \label{t1}
The five two-step nilpotent Lie algebras $N^{8, 2}_i$ for $i= 1, \cdots, 5$ in the last section, are mutually non-isomorphic.
\end{theorem}

This proves the corresponding result of \cite{RZ11}, quantitatively.

Now, we determine which of $N^{8, 2}_i$ for $i = 2, 3, 4$ should be assigned as $T^{8, 2}_{4, i}$ for $i = 2, 3, 4$, respectively.   Since their $H$-msg related sequences are the same, we now follow the lexicographic order of their generator relation sequence first, then follow the lexicographic order of their center related sequence.
Hence, we denote $T^{8, 2}_i$ for $i = 4, 3, 2$, as $N^{8, 2}_{4, i}$ for $i = 2, 3, 4$, respectively.   This totally solves the nomenclature problem for $n = 8$ and $p = 2$.

The case that $n=9$ and $p=2$ can be treated similarly.

\begin{theorem} \label{t2}
The five two-step nilpotent Lie algebras $N^{9, 2}_i$ for $i= 1, \cdots, 5$ in the last section, are mutually non-isomorphic.
\end{theorem}
{\bf Proof}  They all satisfies the condition of Lemma \ref{l1}, i.e., Lemma 2.3 of \cite{RZ11}.  The $H$-msg related sequences of $N^{9, 2}_i$ for $i = 1, 2$ are $(1, 1, 1, 1, 1, 1, 2)$.   The center related sequence for $N^{9, 2}_1$ is $(1, 3)$.  The center related sequence for $N^{9, 2}_2$ is $(2, 2)$.
The $H$-msg related sequences of $N^{9, 2}_i$ for $i= 3, 4$, are  $(1, 1, 1, 1, 2, 2, 2)$.   The generator relation sequence of $N^{9, 2}_3$ is $(2, 5)$.   The generator relation sequence of $N^{9, 2}_4$ is $(3, 4)$.   The $H$-msg related sequence of $N^{9, 2}_5$ is  $(1, 1, 2, 2, 2, 2, 2)$.  The invariants of any two of them are different.   Hence, these five two-step nilpotent Lie algebras are mutually non-isomorphic.
\qed

The generator relation sequences of  $N^{9, 2}_i$ for $i = 1, 2$ are the same as $(2, 2, 3)$.  The center related sequence for $N^{9, 2}_1$ is $(1, 3)$, which precedes $(2, 2)$, the center related sequence for $N^{9, 2}_2$.  Thus we denote $N^{9, 2}_1$ as $T^{9, 2}_{5, 1}$ and $N^{9, 2}_2$ as $T^{9, 2}_{5, 2}$.  The generator relation sequence for $N^{9, 2}_3$ is $(2, 5)$, which precedes $(3, 4)$, the generator relation sequence for $N^{9, 2}_4$.   Thus we denote $N^{9, 2}_3$ as $T^{9, 2}_{4, 1}$ and $N^{9, 2}_4$ as $T^{9, 2}_{4, 2}$.   This completes the nomenclature task for $n = 9$ and $p = 2$.

We now study the case that $n=8$ and $p=3$.

In the last section, we see that for $N^{8, 3}_i$
with $i = 1, 2, 4$, their ranks are all $4$, all the root spaces have dimension $1$, and their $H$-msg related sequences are all the same.  We now need to show that $N^{8, 3}_i$ with $i = 1, 2, 4$ are mutually non-isomorphic, and determine the values of $i$ precisely in $T^{8, 3}_{4, i}$ with $i = 2, 3, 4$ for them.   The generator relation sequence of $N^{8, 3}_1$ is $(2, 3)$, while the generator relation sequence of $N^{8, 3}_2$ and $N^{8, 3}_4$ are $(5)$.   Thus, $N^{8, 3}_1$ is not isomorphic from $N^{8, 3}_2$ and $N^{8, 3}_4$, and we may denote it as $T^{8, 3}_{4, 2}$, as $(2, 5) \prec (5)$.

However, the center related sequences of $N^{8, 3}_2$ and $N^{8, 3}_4$ are the same as $(1, 1, 2)$.   Hence, a new invariant is needed to distinguish them.   Hence, we now define weighted center related sequence for a two-step nilpotent Lie algebra $N$,   Assume that the generating graph $G$ is a $3$-uniform hypergraph, i.e., the Lie bracket operations which define $N$ all have the form $[\x_i, \x_j] = \y_k$.  This assumption is satisfied by $N^{8, 3}_i$ for $i = 1, \cdots, 11$.   For each $i = 1, \cdots, q$, $\x_i$ is a vertex of the generator graph $\hat G$.  Denote the degree of $\x_i$ as a vertex in $\hat G$ as $d(\x_j)$.   For each $k = 1, \cdots, p$, let
$$w(\y_k) = \sum \left\{ d_{\x_i} : (\x_i, \x_j, \y_k) \ \text{is\ a\ hyper-edge\ of\ }G \right\}.$$
We put these numbers in the nondecreasing order as a sequence, and call them the {\bf weighted center related sequence} of $N$.  Under the condition of Lemma 5 of \cite{RZ11}, i.e., Lemma 3 of \cite{RZ17a}, and the $3$-uniform hypergraph assumption, this sequence is also an invariant of $N$, and we may use its lexicographic order to determine the value $i$ in $T^{n, p}_{r, i}$, if the generator relation sequence and center related sequence have already been used.   Now the weighted center related sequence of $N^{8, 3}_2$ is $(4, 4, 6)$, while the weighted center related sequence of $N^{8, 3}_4$ is $(3, 4, 7)$. Thus, $N^{8, 3}_2$ and $N^{8, 3}_4$ are not isomorphic.  Since
$(3, 4, 7) \prec (4, 4, 6)$, we denote  $N^{8, 3}_4$ as  $T^{8, 3}_{4, 3}$, and $N^{8, 3}_2$ as  $T^{8, 3}_{4, 4}$.

Finally, we study the case of $N^{8, 3}_i$ with $i = 6, 9$.    We found that the generator relation sequences, the center related sequences, the weighted center related sequences of these two two-step nilpotent Lie algebras are all the same.   However, the generator graph $\hat G$ of $N^{8, 3}_6$ has a cycle of three, while the generator graph $\hat G$ of $N^{8, 3}_9$ has a cycle of four.   Thus, they are not isomorphic, and $N^{8, 3}_9$ is more complicated in a certain sense.    Thus, we may denote $N^{8, 3}_6$ as $T^{8, 3}_{3, 2}$, and $N^{8, 3}_9$ as $T^{8, 3}_{3, 3}$.  This completes our
nomenclature, and we have the following theorem.

\begin{theorem} \label{t3}
The eleven two-step nilpotent Lie algebras $N^{8, 3}_i$ for $i= 1, \cdots, 11$ in the last section, are mutually non-isomorphic.
\end{theorem}

\section{Identifying Some Indecomposable Two-Step Nilpotent Lie Algebras in Dimension $8$}

In this section, we identify some two-step nilpotent Lie algebras in dimension $8$ being indecomposable.    The tool of generating hypergraphs is used in the proof of Lemma \ref{lem:decomp}.

Given a two-step Lie algebra $N$ with $q$-generators and a center of dimension $p$, and a corresponding basis $(X,I)$ with $X=\{\mathbf x_1,\dots,\mathbf x_q\}$ and $I=\{\mathbf y_1,\dots,\mathbf y_p\}$, we can associate it a third order tensor $\mathcal A(X,I)=(a_{ijk})$ of dimension $q\times q\times p$ with the entries being
\[
a_{ijk}=\alpha_k\ \text{if }[\mathbf x_i,\mathbf x_j]=\sum_{s=1}^p\alpha_s\mathbf y_s
\]
for all $i,j\in\{1,\dots,q\}$ and $k\in\{1,\dots,p\}$. We call it the \textit{representation tensor} of the basis $(X,I)$.

Given a triple of matrices $(A,B,C)\in\mathbb C^{p\times q}\times\mathbb C^{r\times s}\times\mathbb C^{u\times v}$, we have an action of $(A,B,C)$ on a given tensor $\mathcal A\in \mathbb C^{q\times s\times v}$ defined entry-wisely by
\[
\big((A,B,C)\cdot\mathcal A\big)_{ijk}:=\sum_{x=1}^q\sum_{b=1}^s\sum_{c=1}^vA_{ix}B_{jb}C_{kc}a_{xbc}
\]
for all $(i,j,k)\in\{1,\dots,p\}\times\{1,\dots,r\}\times\{1,\dots,u\}$.

Given a Lie algebra $N$ with basis $Z=\{\mathbf z_1,\dots,\mathbf z_n\}$. The \textit{coefficient tensor} $\mathcal A(Z)\in\mathbb C^{n\times n\times n}$ is defined entry-wisely via
\[
[\mathbf z_i,\mathbf z_j]=\sum_{k=1}^n(\mathcal A(Z))_{ijk}\mathbf z_k.
\]
%-----------------------------------------------------------------------------------
\begin{lemma}\label{lem:base-coefficient}
Let $\hat Z=\{\hat{\mathbf z}_1,\dots,\hat{\mathbf z}_n\}$ be another basis of $N$ and $T\in\mathbb C^{n\times n}$ be the basis change matrix from $Z$ to $\hat Z$. Then, we have that
\[
\mathcal A(\hat Z)=\left(T^{-\mathsf{T}},T^{-\mathsf{T}},T\right)\cdot\mathcal A(Z).
\]
\end{lemma}
{\bf Proof}
Let $\mathbf t_i$ be the $i$-th column of $T$. We have
\[
\mathbf z_i = \sum_{r=1}^n(\mathbf t_i)_r\hat{\mathbf z}_r
\]
for all $i\in\{1,\dots,n\}$.

It follows from the bilinearity of the Lie bracket that
\[
\sum_{r,k=1}^n(\mathbf t_i)_r(\mathbf t_j)_k[\hat{\mathbf z}_r,\hat{\mathbf z}_k] = \sum_{k=1}^n(\mathcal A(Z))_{ijk}\mathbf z_k= \sum_{k=1}^n(\mathcal A(Z))_{ijk}\sum_{t=1}^n(\mathbf t_k)_t\hat{\mathbf z}_t.
\]
On the other hand, we have
\[
[\hat{\mathbf z}_r,\hat{\mathbf z}_k]=\sum_{t=1}^n(\mathcal A(\hat Z))_{rkt}\hat{\mathbf z}_t.
\]
As a result, we have for all possible $(i,j,t)$
\[
\sum_{r,k=1}^n(\mathbf t_i)_r(\mathbf t_j)_k (\mathcal A(\hat Z))_{rkt}= \sum_{k=1}^n(\mathcal A(Z))_{ijk} (\mathbf t_k)_t.
\]
Thus, with $E$ being the identity matrix of appropriate dimension,
\[
(T^\mathsf{T},T^\mathsf{T},E)\cdot\mathcal A(\hat Z)=(E,E,T)\cdot\mathcal A(Z).
\]
By the nonsingularity of the matrix $T$ and the associativity of the action on the tensors defined above, we have that
\begin{align*}
\mathcal A(\hat Z)&=(T^{-\mathsf{T}},T^{-\mathsf{T}},E)\cdot\big(
(T^\mathsf{T},T^\mathsf{T},E)\cdot\mathcal A(\hat Z)\big)\\
&=(T^{-\mathsf{T}},T^{-\mathsf{T}},E)\cdot\big((E,E,T)\cdot\mathcal A(Z)\big)\\
&=(T^{-\mathsf{T}},T^{-\mathsf{T}},T)\cdot\mathcal A(Z).
\end{align*}
The conclusion then follows.
\qed
%-----------------------------------------------------------------------------------
\begin{lemma}\label{lem:base}
Let $(\hat X=\{\hat{\mathbf x}_1,\dots,\hat{\mathbf x}_q\},\hat I=\{\hat{\mathbf y}_1,\dots,\hat{\mathbf y}_p\})$ be another basis of the two step Lie algebra $N$ and
\begin{equation}\label{eq:basis}
Q=\begin{bmatrix}
S&0\\ P& C
\end{bmatrix}
\end{equation}
with $S\in\mathbb C^{q\times q}$ and $C\in\mathbb C^{p\times p}$ be the basis change matrix from $(X,I)$ to $(\hat X,\hat I)$. Then, we have that
\[
\mathcal A(\hat X,\hat I)=\left(S^{-\mathsf{T}},S^{-\mathsf{T}},C\right)\cdot\mathcal A(X,I).
\]
\end{lemma}
{\bf Proof}
It follows from Lemma~\ref{lem:base-coefficient} and the fact that
\[
Q^{-1}=\begin{bmatrix}S^{-1}&0\\ M&C^{-1} \end{bmatrix}
\]
for some matrix $M$.
\qed

Note that both the matrices $S$ and $C$ in $Q$ are nonsingular. On the other hand, any pair of nonsingular matrices $S$ and $C$ will compose a basis change matrix $Q$. 

A third order tensor $\mathcal A\in \mathbb C^{q\times q\times p}$ is in \textit{block diagonal format $(S,T)$} with $S\subseteq\{1,\dots,q\}$ and $T\subseteq\{1,\dots,p\}$ if the only possible nonzero entries of $\mathcal A$ occur at
\[
(i,j,k)\in S\times S\times T\cup S^\complement\times S^\complement \times T^\complement.
\]

%-------------------------------------
\begin{lemma}\label{lem:decomp}
Let all notation be as above. Then the Lie algebra $N$ is decomposable if and only if there is a basis $(X,I)$ such that the representation tensor $\mathcal A(X,I)$ is in block diagonal format $(S,T)$ for some nonempty proper subset $S\subset\{1,\dots,q\}$.
\end{lemma}
{\bf Proof}
In this case, the generating hypergraph of $N$ has two disconnected parts, induced by the vertex sets $S\cup T$ and $S^\complement\cup T^\complement$. Thus, by Proposition~\ref{p3}, the conclusion follows.
\qed

The decomposable case with a nonempty proper subset $S\subset\{1,\dots,q\}$ and $T=\{1,\dots,p\}$ can be identified easily. %Actually, this case gives a Lie algebra being a direct sum of two nontrivial Lie algebras.

Let $A_i:=(a_{\cdot\cdot i})\in\mathbb C^{q\times q}$ be the $i$-th slice matrix of the tensor $\mathcal A(X,I)$ for all $i\in\{1,\dots,q\}$. The rank of the matrix $[A_1,\dots,A_p]\in\mathbb C^{q\times pq}$ is called the \textit{marginal rank} of $N$. Since the representation tensor is changed along the bases change by the rule in Lemma~\ref{lem:base}, obviously, the marginal rank is an intrinsic quantity which is independent of bases.
%-------------------------------------
\begin{lemma}\label{lem:marginal}
Let all notation be as above. There is a basis $(X,I)$ such that the representation tensor $\mathcal A(X,I)$ is in block diagonal format $(S,T)$ for some nonempty proper subset $S\subset\{1,\dots,q\}$ and $T=\{1,\dots,p\}$ if and only if the marginal rank of $N$ is strictly smaller than $q$. In this case, the Lie algebra $N$ is a direct sum of $N^{s+p,p}$ and $q-s$ copies of $N^1$, and hence decomposable.
\end{lemma}
{\bf Proof}
For the necessity, under this basis, the matrix $[A_1,\dots,A_p]$ has only nonzero entries in the rows indexed by the set $S$. Thus, the marginal rank is obviously not greater than $|S|<q$.

For the sufficiency, if the marginal rank is not greater than a positive number $s<q$, there is an invertible matrix $Q$ such that $Q[A_1,\dots,A_p]$ has only nonzero entries in the rows indexed by a proper nonempty subset $S$ of $\{1,\dots,q\}$ with cardinality $|S|\leq s$. Also note that each $A_i$ is skew-symmetric by the skew-symmetry of the Lie bracket. Therefore, $QA_iQ^\mathsf{T}$ is also skew-symmetric. Since $QA_i$ has zero rows indexed by $S^\complement$, we have that $QA_iQ^\mathsf{T}$ also has zero columns indexed by $S^\complement$. On the other hand, the $i$-th slice of the tensor $(Q,Q,E)\cdot \mathcal A(X,I)$ is exactly $QA_iQ^\mathsf{T}$ for all $i\in\{1,\dots,p\}$. Consequently, the tensor $(Q,Q,E)\cdot \mathcal A(X,I)$ has the desired block diagonal format.
\qed

Block diagonal formats as in Lemma~\ref{lem:marginal} are called \textit{trivial}. The block diagonal formats with nonempty proper subsets $T\subset\{1,\dots,p\}$ are \textit{nontrivial}.

\begin{lemma}\label{lem:rank}
Let $p=2$ and the marginal rank of $N$ be $q$. Let $\mathcal A_{\cdot\cdot i}\in\mathbb C^{q\times q}$ be the matrix slice of $\mathcal A$ by fixing the third index to $i$ for $i=1,2$. Then, if
\begin{equation}\label{eq:rank}
\operatorname{rank}(\mathcal A_{\cdot\cdot 1})+\operatorname{rank}(\mathcal A_{\cdot\cdot 2})>q,
\end{equation}
then $\mathcal A$ cannot be block diagonal format $(S,T)$ for any nonempty proper subset $S\subset\{1,\dots,q\}$. Therefore, if \eqref{eq:rank} is satisfied for all the representation tensor of the Lie algebra $N$, then $N$ is indecomposable.
\end{lemma}
{\bf Proof}
By Lemma~\ref{lem:marginal}, only nontrivial block diagonal format with nonempty proper subset $T\subset\{1,2\}$ is possible. The rest follows from the definitions.
\qed

\begin{lemma}\label{lem-two-step}
Let $p=2$ and the marginal rank of $N$ be $q$, $A=\mathcal A_{\cdot\cdot 1}$ and $B=\mathcal A_{\cdot\cdot 2}$ for a representation tensor. If
\begin{equation}\label{eq:rank-2}
\operatorname{rank}(c_{11}A+c_{12}B)+\operatorname{rank}(c_{21}A+c_{22}B)>q
\end{equation}
for all nonsingular matrix $C\in\mathbb C^{2\times 2}$, then the Lie algebra $N$ is indecomposable.
\end{lemma}
{\bf Proof}
Let $(\hat X,\hat I)$ be another basis of $N$.
Let $S\in\mathbb C^{q\times q}$ be the basis change matrix from $X$ to $\hat X$, and $C$ be that from $I$ to $\hat I$. It follows from Lemma~\ref{lem:base} that
\[
\mathcal A(\hat X,\hat I)=(S^{-\mathsf{T}},S^{-\mathsf{T}},C)\cdot\mathcal A(X,I).
\]
Let $\hat A, \hat B$ be the slice matrices of $\mathcal A(\hat X,\hat I)$. We have that
\[
\hat A=S^\mathsf{-T}(c_{11}A+c_{12}B)S^{-1}\ \text{and }\hat B=S^\mathsf{-T}(c_{21}A+c_{22}B)S^{-1}.
\]
Since $S$ is nonsingular, we have that
\[
\operatorname{rank}(\hat A)=\operatorname{rank}(c_{11}A+c_{12}B)\ \text{and }\operatorname{rank}(\hat B)=\operatorname{rank}(c_{21}A+c_{22}B).
\]
The result then follows from Lemma~\ref{lem:rank}.
\qed

\begin{theorem}\label{prop-n82}
All the Lie algebras $N^{8,2}_i$ for $i=1,\dots,5$ are indecomposable.
\end{theorem}
{\bf Proof}
The matrices $A$ and $B$ for $N^{8,2}_1$ are
\[
A=\begin{bmatrix}0&1&0&0&0&0\\ -1&0&0&0&0&0\\ 0&0&0&0&0&0\\ 0&0&0&0&0&0\\0&0&0&0&0&1\\ 0&0&0&0&-1&0\end{bmatrix}\ \text{and }B=\begin{bmatrix}0&0&0&0&0&0\\ 0&0&0&0&0&0\\ 0&0&0&1&0&0\\ 0&0&-1&0&0&0\\0&0&0&0&0&1\\ 0&0&0&0&-1&0\end{bmatrix}
\]
It is  easy to see that the marginal rank is $q=6$. Then by Lemma~\ref{lem:marginal}, it cannot be written as a direct sum of $N^{s,2}$ and $8-s$ copies of $N^1$. In the following, we consider nontrivial block diagonal formats.
Obviously, $\operatorname{rank}(A)=4$ and $\operatorname{rank}(B)=4$. Moreover,
\[
\operatorname{rank}(A+\gamma B)\geq 6\ \text{and }\operatorname{rank}(\tau A+ B)\geq 6
\]
for any nonzero $\gamma$ and $\tau$. By Lemma~\ref{lem-two-step}, $N^{8,2}_1$ is indecomposable.

The matrices for $N^{8,2}_2$ are
\[
A=\begin{bmatrix}0&0&0&0&0&-1\\ 0&0&0&0&-1&0\\ 0&0&0&0&0&0\\ 0&0&0&0&0&0\\0&1&0&0&0&0\\ 1&0&0&0&0&0\end{bmatrix}\ \text{and }B=\begin{bmatrix}0&0&0&0&0&0\\ 0&0&0&0&0&0\\ 0&0&0&0&-1&0\\ 0&0&0&0&0&-1\\0&0&1&0&0&0\\ 0&0&0&1&0&0\end{bmatrix}.
\]

The matrices for $N^{8,2}_3$ are
\[
A=\begin{bmatrix}0&1&0&0&0&0\\ -1&0&0&0&0&0\\ 0&0&0&0&0&0\\ 0&0&0&0&0&0\\0&0&0&0&0&-1\\ 0&0&0&0&1&0\end{bmatrix}\ \text{and }B=\begin{bmatrix}0&0&0&0&0&0\\ 0&0&0&0&0&0\\ 0&0&0&0&0&1\\ 0&0&0&0&-1&0\\0&0&0&1&0&0\\ 0&0&-1&0&0&0\end{bmatrix}.
\]

The matrices for $N^{8,2}_4$ are
\[
A=\begin{bmatrix}0&1&0&0&0&0\\ -1&0&0&0&0&0\\ 0&0&0&0&0&1\\ 0&0&0&0&-1&0\\0&0&0&1&0&0\\ 0&0&-1&0&0&0\end{bmatrix}\ \text{and }B=\begin{bmatrix}0&0&0&0&0&0\\ 0&0&0&0&0&0\\ 0&0&0&0&0&0\\ 0&0&0&0&0&0\\0&0&0&0&0&-1\\ 0&0&0&0&1&0\end{bmatrix}.
\]

The matrices for $N^{8,2}_5$ are
\[
A=\begin{bmatrix}0&0&0&0&0&1\\ 0&0&0&0&-1&0\\ 0&0&0&1&0&0\\ 0&0&-1&0&0&0\\0&1&0&0&0&0\\ -1&0&0&0&0&0\end{bmatrix}\ \text{and }B=\begin{bmatrix}0&0&0&0&0&0\\ 0&0&0&0&0&0\\ 0&0&0&0&0&-1\\ 0&0&0&0&1&0\\0&0&0&-1&0&0\\ 0&0&1&0&0&0\end{bmatrix}.
\]

Similar arguments as $N^{8,2}_1$ will show that $N^{8,2}_i$ for all $i=2,\dots,5$ are all indecomposable.
\qed

\section{Identifying Some Decomposable Two-Step Nilpotent Lie Algebras in Dimension $9$}

For $n=9$ and $p=4$, Table 2 of Galitski and Timashev \cite{GT99} listed $35$ non-isomorphic two-step nilpotent Lie algebras.   The No. 82 algebra is defined by
$$[\x_1, \x_2] = \y_2, [\x_1, \x_3] = \y_3, [\x_2, \x_3] = \y_4, [\x_4, \x_5] = \y_1.$$
The vertex set of the generating hypergraph $G$ of this two-step nilpotent Lie algebra can be partitioned into two parts $\{ \x_1, \x_2, \x_3, \y_2, \y_3, \y_4 \}$ and $\{ \x_4, \x_5, \y_1 \}$, which are not connected.    By Proposition \ref{p3}, this two-step nilpotent Lie algebra is decomposable.

Galitski and Timashev \cite{GT99} indicated that for $n=9$ and $p=4$, there are seven families of non-isomorphic two-step nilpotent Lie algebras.     Each of the first three families depends on two parameters.   The fourth, fifth and sixth families contain $6$, $6$ and $15$ algebras respectively.   In Table 8 of \cite{GT99}, $44$ non-isomorphic two-step nilpotent Lie algebras of $n=9$ and $p=4$ were listed.    The No. 44 algebra of Table 8 of \cite{GT99}, is defined by
$$[\x_1, \x_2] = \y_3, [\x_1, \x_5] = \y_1, [\x_2, \x_6] = \y_1, [\x_3, \x_4] = \y_2.$$
The vertex set of the generating hypergraph $G$ of this two-step nilpotent Lie algebra can be partitioned into two parts $\{ \x_1, \x_2, \x_5, \x_6, \y_1, \y_3 \}$ and $\{ \x_3, \x_4, \y_2 \}$, which are not connected.    By Proposition \ref{p3}, this two-step nilpotent Lie algebra is decomposable.

It will be a challenging work to identify which of the remaining two-step nilpotent Lie algebras in these two cases are indecomposable.

\bigskip

{\bf Acknowledgment}   We are thankful to Professors Chengming Bai and Shaoqiang Deng for their encouragements.   We are especially thankful to Dr. Zaili Yan who explained various issues related with this research topic in details.

%%%%%%%%%%%%%%%%%%%%%%%%%%%%%%%

\end{document}